\documentclass[a4, 12pt]{amsart}
\usepackage{amssymb}
\usepackage{amstext}
\usepackage{amsmath}
\usepackage{amscd}
\usepackage{amsthm}
\usepackage{amsfonts}
\usepackage{enumerate}
\usepackage{graphicx}
\usepackage{latexsym}

\theoremstyle{plain}
\newtheorem{thm}{Theorem}[section]
\newtheorem*{thm*}{Theorem}
\newtheorem*{cor*}{Corollary}
\newtheorem*{prop*}{Proposition}
\newtheorem*{thma}{Theorem A}
\newtheorem*{thmb}{Theorem B}

\newtheorem{prop}[thm]{Proposition}
\newtheorem{lem}[thm]{Lemma}

\newtheorem{cor}[thm]{Corollary}

\newtheorem*{claim*}{Claim}

\theoremstyle{definition}
\newtheorem{dfn}[thm]{Definition}
\newtheorem{rem}[thm]{Remark}
\newtheorem{ex}[thm]{Example}

\newtheorem*{conj*}{Conjecture}

\newtheorem*{conv}{{\sc Convention}}
\newtheorem*{ac}{{\sc Acknowledgments}}

\theoremstyle{remark}

\numberwithin{equation}{thm}
\def\Hom{\operatorname{Hom}}

\def\Ext{\operatorname{Ext}}
\def\Tor{\operatorname{Tor}}

\def\mod{\operatorname{mod}}
\def\res{\operatorname{res}}

\def\add{\operatorname{add}}
\def\ind{\operatorname{ind}}
\def\CM{\operatorname{CM}}
\def\proj{\operatorname{proj}}

\def\NF{\operatorname{NF}}

\def\level{\operatorname{step}}

\def\m{\mathfrak m}

\def\p{\mathfrak p}
\def\q{\mathfrak q}
\def\r{\mathfrak r}

\def\Z{\Bbb Z}
\def\C{\Bbb C}

\def\Supp{\operatorname{Supp}}

\def\Ann{\operatorname{Ann}}

\def\height{\operatorname{ht}}
\def\grade{\operatorname{grade}}

\def\Spec{\operatorname{Spec}}
\def\Sing{\operatorname{Sing}}

\def\G{{\mathcal G}}

\def\X{{\mathcal X}}
\def\Y{{\mathcal Y}}

\begin{document}

\setlength{\baselineskip}{15pt}

\title[Modules in resolving subcategories]{Modules in resolving subcategories which are free on the punctured spectrum}
\author{Ryo Takahashi}
\address{Department of Mathematical Sciences, Faculty of Science, Shinshu University, 3-1-1 Asahi, Matsumoto, Nagano 390-8621, Japan}
\email{takahasi@math.shinshu-u.ac.jp}
\keywords{resolving subcategory, resolving closure, nonfree locus, Cohen-Macaulay ring, maximal Cohen-Macaulay module, countable Cohen-Macaulay representation type, totally reflexive module}
\subjclass[2000]{13C05, 16D90, 16G60, 16G50, 13C14, 13C13}
\begin{abstract}
Let $R$ be a commutative noetherian local ring, and let $\X$ be a resolving subcategory of the category of finitely generated $R$-modules.
In this paper, we study modules in $\X$ by relating them to modules in $\X$ which are free on the punctured spectrum of $R$.
We do this by investigating nonfree loci and establishing an analogue of the notion of a level in a triangulated category which has been introduced by Avramov, Buchweitz, Iyengar and Miller.
As an application, we prove a result on the dimension of the nonfree locus of a resolving subcategory having only countably many nonisomorphic indecomposable modules in it, which is a generalization of a theorem of Huneke and Leuschke.
\end{abstract}
\maketitle
\section{Introduction}

In the 1960s, Auslander and Bridger \cite{AB} introduced the notion of a resolving subcategory of an abelian category with enough projectives.
They proved that in the category of finitely generated modules over a left and right noetherian ring, the full subcategory consisting of all modules of Gorenstein dimension zero, which are now also called totally reflexive modules, is resolving.

Let $R$ be a commutative noetherian ring, and let $\mod R$ denote the category of finitely generated $R$-modules.
A lot of important full subcategories of $\mod R$ are known to be resolving.
As trivial examples, $\mod R$ itself and the full subcategory $\proj R$ of $\mod R$ consisting of all projective modules are resolving.
If $R$ is a Cohen-Macaulay local ring, then the full subcategory $\CM(R)$ of $\mod R$ consisting of all maximal Cohen-Macaulay $R$-modules is resolving; see Example \ref{ex} for details and other examples of a resolving subcategory.

Let $R$ be a local ring, and let $\X$ be a resolving subcategory of $\mod R$.
In the present paper, we study modules in $\X$ by relating them to modules in $\X$ which are free on the punctured spectrum of $R$.
A key role is played by the nonfree loci of $R$-modules and subcategories of $\mod R$, which are certain closed and specialization-closed subsets of $\Spec R$, respectively.

To be more precise, from each module $X\in\X$ we construct another module $X'\in\X$ which is free on the punctured spectrum of $R$, and count the (minimum) number of steps required to construct $X'$ from $X$.
We denote the number by $\level(X,X')$.
(The precise definition will be given in Definition \ref{deflevel}.)
This invariant is an analogue of a level in a triangulated category which has been introduced by Avramov, Buchweitz, Iyengar and Miller \cite{ABIM}.

We denote by $\NF(X)$ and $\NF(\X)$ the nonfree loci of an $R$-module $X$ and a subcategory $\X$ of $\mod R$, respectively (cf. Definition \ref{defnf}).
The main result of this paper is the following, which will be proved in Corollary \ref{locthm2}.

\begin{thma}
Let $R$ be a commutative noetherian local ring, and let $\X$ be a resolving subcategory of $\mod R$.
Then for every nonfree $R$-module $X\in\X$, there exists a nonfree $R$-module $X'\in\X$ satisfying the following two conditions:
\begin{enumerate}[\rm (1)]
\item
$\level(X,X')\le 2\dim\NF(X)$, and
\item
$X'$ is free on the punctured spectrum of $R$.
\end{enumerate}
\end{thma}

As an application, we consider how many nonisomorphic indecomposable modules are in $\X$.
We will prove the following result in Corollary \ref{corthm3}.

\begin{thmb}
Let $R$ be a commutative noetherian local ring which is either complete or has uncountable residue field.
Let $\X$ be a resolving subcategory of $\mod R$ in which there are only countably many nonisomorphic indecomposable $R$-modules.
Then $\dim\NF(\X)\le 1$.
\end{thmb}

This theorem recovers a theorem of Huneke and Leuschke \cite{HL2} which proves a conjecture of Schreyer \cite{S} in the 1980s.

\begin{conv}
Throughout this paper, let $R$ be a commutative noetherian ring.
All $R$-modules considered in this paper are assumed to be finitely generated.
We denote by $\mod R$ the category of finitely generated $R$-modules.
By a {\em subcategory} of $\mod R$, we always mean a full subcategory of $\mod R$ which is closed under isomorphisms.
We freely use basic definitions and results in commutative algebra which are stated in \cite{BH}.
\end{conv}

\begin{ac}
The author is grateful to Yuji Kamoi, Takesi Kawasaki and Koji Nishida for valuable discussions and helpful advices.
He also thanks the referee for useful comments and suggestions.
He was partially supported by Grant-in-Aid for Young Scientists (B) 19740008 from JSPS. 
\end{ac}

\section{Foundations}

In this section, we define the resolving closures and the nonfree loci of an $R$-module and a subcategory of $\mod R$, and study their basic properties.
We begin with recalling the definition of a syzygy.

\begin{dfn}
\begin{enumerate}[(1)]
\item
Let $n$ be a nonnegative integer, and let $M$ be an $R$-module.
If there exists an exact sequence
$$
0 \to N \to P_{n-1} \to \cdots \to P_0 \to M \to 0
$$
of $R$-modules where $P_i$ is a projective $R$-module for every $0\le i\le n-1$, then we call $N$ the $n$th {\em syzygy} of $M$, and denote it by $\Omega^nM$.
Note that the $n$th syzygy of a given $R$-module is not uniquely determined; it is uniquely determined up to projective summand.
We simply write $\Omega^1M=\Omega M$.
\item
In the case where $R$ is local, the $R$-module $M$ admits a minimal free resolution
$$
\cdots \overset{\partial_{n+1}}{\to} F_n \overset{\partial_n}{\to} \cdots \overset{\partial_1}{\to} F_0 \overset{\partial_0}{\to} M \to 0.
$$
Then we define the $n$th {\em syzygy} of $M$ as the image of $\partial_n$, and denote it by $\Omega^nM$.
The $n$th syzygy of a given $R$-module is uniquely determined up to isomorphism since so is a minimal free resolution.
Whenever $R$ is local, we define the $n$th syzygy of an $R$-module by using its minimal free resolution.
\end{enumerate}
\end{dfn}

Next let us recall the definition of a resolving subcategory.

\begin{dfn}
A subcategory $\X$ of $\mod R$ is called {\em resolving} if $\X$ satisfies the following conditions.
\begin{enumerate}[(1)]
\item
$\X$ contains all projective $R$-modules.
\item
$\X$ is closed under direct summands: if $M$ is in $\X$ and $N$ is a direct summand of $M$, then $N$ is also in $\X$.
\item
$\X$ is closed under extensions: for any exact sequence $0\to L\to M\to N\to 0$ in $\mod R$, if $L$ and $N$ are in $\X$, then so is $M$.
\item
$\X$ is closed under kernels of epimorphisms: for any exact sequence $0\to L\to M\to N\to 0$ in $\mod R$, if $M$ and $N$ are in $\X$, then so is $L$.
\end{enumerate}
\end{dfn}

A resolving subcategory is a subcategory such that any two ``minimal'' resolutions of a module by modules in it have the same length; see \cite[Lemma (3.12)]{AB}.

The closedness under kernels of epimorphisms can be replaced with a weaker condition of the closedness under syzygies.

\begin{rem}\cite[Lemma 3.2]{Y2}
A subcategory $\X$ of $\mod R$ is resolving if and only if $\X$ satisfies the following conditions.
\begin{enumerate}[(1)]
\item
$\X$ contains all projective $R$-modules.
\item
$\X$ is closed under direct summands.
\item
$\X$ is closed under extensions.
\item
$\X$ is closed under syzygies: if $M$ is in $\X$, then so is $\Omega M$.
\end{enumerate}
\end{rem}

A lot of important subcategories of $\mod R$ are known to be resolving.
Here, let us make a list of examples.

\begin{ex}\label{ex}
\begin{enumerate}[(1)]
\item
It is trivial that the subcategory $\mod R$ of $\mod R$ is resolving.
\item\label{proj}
It is obvious that the subcategory $\proj R$ of $\mod R$ consisting of all projective $R$-modules is resolving.
\item\label{grade}
Let $I$ be an ideal of $R$.
Then the subcategory of $\mod R$ consisting of all $R$-modules $M$ with $\grade(I,M)\ge\grade(I,R)$ is resolving.
This can be shown by using the equality $\grade(I,M)=\inf\{i\in\Z\mid\Ext_R^i(R/I,M)\ne 0\}$.
\item\label{cm}
Let $R$ be a Cohen-Macaulay local ring.
Then, letting $I$ be the maximal ideal of $R$ in (\ref{grade}), we see that the subcategory $\CM(R)$ of $\mod R$ consisting of all maximal Cohen-Macaulay $R$-modules is resolving.
\item\label{tcr}
An $R$-module $C$ is called {\em semidualizing} if the natural homomorphism $R\to\Hom_R(C,C)$ is an isomorphism and $\Ext_R^i(C,C)=0$ for every $i>0$.
An $R$-module $M$ is called {\em totally $C$-reflexive}, where $C$ is a semidualizing $R$-module, if the natural homomorphism $M\to\Hom_R(\Hom_R(M,C),C)$ is an isomorphism and $\Ext_R^i(M,C)=\Ext_R^i(\Hom_R(M,C),C)=0$ for every $i>0$.
The subcategory $\G_C(R)$ of $\mod R$ consisting of all totally $C$-reflexive $R$-modules is resolving by \cite[Theorem 2.1]{ATY}.
\item\label{tr}
A totally $R$-reflexive $R$-module is simply called {\em totally reflexive}.
The subcategory $\G(R)$ of $\mod R$ consisting of all totally reflexive $R$-modules is resolving by (\ref{tcr}); see also \cite[(3.11)]{AB}.
\item
Let $n$ be a nonnegative integer, and let $K$ be an $R$-module (which is not necessarily finitely generated).
Then the subcategory of $\mod R$ consisting of all $R$-modules $M$ with $\Tor_i^R(M,K)=0$ for $i>n$ (respectively, $i\gg 0$) and the subcategory of $\mod R$ consisting of all $R$-modules $M$ with $\Ext_R^i(M,K)=0$ for $i>n$ (respectively, $i\gg 0$) are both resolving.
\item
Let $R$ be a local ring.
We say that an $R$-module $M$ is {\em bounded} if there is an integer $s$ such that $\beta_i^R(M)\le s$ for all $i\ge 0$, where $\beta_i^R(M)$ denotes the $i$th Betti number of $M$.
The subcategory of $\mod R$ consisting of all bounded $R$-modules is resolving.
This can be shown by using the equality $\beta_i^R(M)=\dim_k\Tor_i^R(M,k)$, where $k$ is the residue field of $R$.
\item\label{cx}
Let $R$ be local.
We say that an $R$-module $M$ has {\em complexity} $c$ if $c$ is the least nonnegative integer $d$ such that there exists a real number $r$ satisfying the inequality $\beta_i^R(M)\le ri^{d-1}$ for $i\gg 0$.
The subcategory of $\mod R$ consisting of all $R$-modules having finite complexity is resolving by \cite[Proposition 4.2.4]{A2}.
\item
Let $R$ be local.
We say that an $R$-module $M$ has {\em lower complete intersection zero} if $M$ is totally reflexive and has finite complexity.
The subcategory of $\mod R$ consisting of all $R$-modules of lower complete intersection dimension zero is resolving by (\ref{tr}) and (\ref{cx}); see also \cite[Lemma 6.3.1]{A}.
\end{enumerate}
\end{ex}

Now we define the resolving closures of a subcategory of $\mod R$ and an $R$-module.

\begin{dfn}
For a subcategory $\X$ of $\mod R$, we denote by $\res\X$ (or $\res_R\X$ when there is some fear of confusion) the resolving sucategory of $\mod R$ generated by $\X$, namely, the smallest resolving subcategory of $\mod R$ containing $\X$.
If $\X$ consists of a single module $X$, then we simply write $\res X$ (or $\res_RX$).
\end{dfn}

\begin{rem}
\begin{enumerate}[(1)]
\item
Let $\{\X_\lambda\}_{\lambda\in\Lambda}$ be a family of resolving subcategories of $\mod R$.
Then the intersection $\bigcap_{\lambda\in\Lambda}\X_\lambda$ is also a resolving subcategory of $\mod R$.
Therefore, for every subcategory $\X$ of $\mod R$, the smallest resolving subcategory of $\mod R$ containing $\X$ exists.
\item
Let $\X,\Y$ be subcategories of $\mod R$.
If $\X\subseteq\Y$, then $\res\X\subseteq\res\Y$.
\item
A subcategory $\X$ of $\mod R$ is resolving if and only if $\X=\res\X$.
In particular, $\res\X=\res(\res\X)$ for every subcategory $\X$ of $\mod R$.
\end{enumerate}
\end{rem}

Next we recall the definition of the nonfree locus of an $R$-module and define the nonfree locus of a subcategory of $\mod R$.

\begin{dfn}\label{defnf}
\begin{enumerate}[(1)]
\item
We denote by $\NF(X)$ (or $\NF_R(X)$) the {\em nonfree locus} of an $R$-module $X$, namely, the set of prime ideals $\p$ of $R$ such that the $R_\p$-module $X_\p$ is nonfree.
\item
We define the {\em nonfree locus} of a subcategory $\X$ of $\mod R$ as the union of $\NF(X)$ where $X$ runs through all (nonisomorphic) $R$-modules in $\X$, and denote it by $\NF(\X)$ (or $\NF_R(\X)$).
\end{enumerate}
\end{dfn}

\begin{rem}
\begin{enumerate}[(1)]
\item
For a subcategory $\X$ of $\mod R$, one has $\NF(\X)=\emptyset$ if and only if $\X$ is contained in $\proj R$.
In particular, one has $\NF(X)=\emptyset$ for an $R$-module $X$ if and only if $X$ is projective.
\item
Let $R$ be a local ring with maximal ideal $\m$.
Then an $R$-module $X$ is nonfree if and only if $\m$ is in $\NF(X)$.
\item
Let $\X,\Y$ be subcategories of $\mod R$.
If $\X\subseteq\Y$, then $\NF(\X)\subseteq\NF(\Y)$.
\end{enumerate}
\end{rem}

\begin{ex}\label{cmex}
Let $R$ be a Cohen-Macaulay local ring.
Then the nonfree locus $\NF(\CM(R))$ coincides with the singular locus $\Sing R$ of $R$.

In fact, for every prime ideal $\p$ in $\NF(\CM(R))$ there exists a maximal Cohen-Macaulay $R$-module $X$ such that the $R_\p$-module $X_\p$ is not free.
Hence $X_\p$ is a nonfree maximal Cohen-Macaulay $R_\p$-module, which implies that the local ring $R_\p$ is singular.
On the other hand, each prime ideal $\p$ in $\Sing R$ belongs to the nonfree locus of the maximal Cohen-Macaulay $R$-module $\Omega^d(R/\p)$, where $d=\dim R$.
\end{ex}

The nonfree locus of a module can be described as the support of an Ext module.

\begin{prop}\label{nfsupp}
Let $\sigma:0 \to Y \to P \to X \to 0$ be an exact sequence of $R$-modules such that $P$ is projective.
Then one has $\NF(X)=\Supp\Ext_R^1(X,Y)$.
Hence $\NF(X)=\Supp\Ext_R^1(X,\Omega X)$.
\end{prop}

\begin{proof}
For a prime ideal $\p$ in $\Supp\Ext_R^1(X,Y)$, the module $\Ext_{R_\p}^1(X_\p,Y_\p)$ is nonzero.
In particular, $X_\p$ is a nonfree $R_\p$-module, and hence $\p$ is in $\NF(X)$.
Conversely, let $\p$ be a prime ideal in $\NF(X)$.
Localizing $\sigma$ at $\p$, we obtain an exact sequence $\sigma_\p:0 \to Y_\p \to P_\p \to X_\p \to 0$ of $R_\p$-modules.
Since $X_\p$ is not free, this exact sequence $\sigma_\p$ does not split, hence this defines a nonzero element of the module $\Ext_{R_\p}^1(X_\p,Y_\p)$.
Thus $\Ext_R^1(X,Y)_\p$ is nonzero, that is, $\p$ is in $\Supp\Ext_R^1(X,Y)$.
\end{proof}

Recall that a subset $Z$ of $\Spec R$ is called {\em specialization-closed} provided that if $\p\in Z$ and $\q\in\Spec R$ with $\p\subseteq\q$ then $\q\in Z$.
Note that every closed subset of $\Spec R$ is specialization-closed.

\begin{cor}\label{cl-spcl}
\begin{enumerate}[\rm (1)]
\item
The nonfree locus of an $R$-module is a closed subset of $\Spec R$ in the Zariski topology.
\item
The nonfree locus of a subcategory of $\mod R$ is specialization-closed.
\end{enumerate}
\end{cor}

\begin{proof}
(1) It is seen from Proposition \ref{nfsupp} that $\NF(X)=\Supp\Ext_R^1(X,\Omega X)$ for an $R$-module $X$.
As $\Ext_R^1(X,\Omega X)$ is a finitely generated $R$-module, the subset $\NF(X)$ of $\Spec R$ is closed.

(2) It is easy to see that in general any union of closed subsets of $\Spec R$ is specialization-closed.
Hence this statement follows from (1).
\end{proof}

Let $Z$ be a closed subset of $\Spec R$.
Then $Z=V(I)$ for some ideal $I$ of $R$.
We call such an ideal $I$ the {\em defining ideal} of $Z$.
This is uniquely determined up to radical.

\section{Inductive construction of resolving closures}

In this section, we build a filtration of subcategories in the resolving closure of a subcategory of $\mod R$, and inductively construct the resolving closure.
This is an imitation of the notion of thickenings in the thick closure of a subcategory of a triangulated category, which were introduced by Avramov, Buchweitz, Iyengar and Miller \cite{ABIM}.
Using this construction of a resolving closure, we can obtain several properties of a resolving closure and its nonfree locus.

The {\em additive closure} $\add\X$ (or $\add_R\X$) of a subcategory $\X$ of $\mod R$ is defined to be the subcategory of $\mod R$ consisting of all direct summands of finite direct sums of modules in $\X$.
Note that $\add\X$ is closed under direct summands and finite direct sums, namely, $R$-modules $M$ and $N$ both belong to $\add\X$ if and only if so does $M\oplus N$.

\begin{dfn}
Let $\X$ be a subcategory of $\mod R$.
For a nonnegative integer $n$, we inductively define a subcategory $\res^n\X$ (or $\res_R^n\X$) of $\mod R$ as follows:
\begin{enumerate}[(1)]
\item
Set $\res^0\X=\add(\X\cup\{ R\})$.
\item
For $n\ge 1$, let $\res^n\X$ be the additive closure of the subcategory of $\mod R$ consisting of all $R$-modules $Y$ having an exact sequence of either of the following two forms:
\begin{align*}
& 0 \to A \to Y \to B \to 0,\\
& 0 \to Y \to A \to B \to 0
\end{align*}
where $A,B\in\res^{n-1}\X$.
\end{enumerate}
If $\X$ consists of a single module $X$, then we simply write $\res^nX$ instead of $\res^n\X$.
\end{dfn}

\begin{rem}
Let $\X,\Y$ be subcategories of $\mod R$, and let $n$ be a nonnegative integer.
Then the following hold.
\begin{enumerate}[\rm (1)]
\item
If $\X\subseteq\Y$, then $\res^n\X\subseteq\res^n\Y$.
\item
One has equalities $\res^n(\add\X)=\res^n\X=\add(\res^n\X)$.
\item
There is an ascending chain $\{0\}\subseteq\res^0\X\subseteq\res^1\X\subseteq\cdots\subseteq\res^n\X\subseteq\cdots\subseteq\res\X$ of subcategories of $\mod R$.
\item
The equality $\res\X=\bigcup_{n\ge 0}\res^n\X$ holds.
\end{enumerate}

The first and second statements follow by definition and induction on $n$.
As to the third statement, since $\res^n\X$ is closed under direct summands, it contains the zero module $0$.
For an $R$-module $M$ in $\res^n\X$ there exists a short exact sequence $0 \to M \overset{=}{\to} M \to 0 \to 0$, which shows that $M$ is in $\res^{n+1}\X$ by definition.
As for the fourth statement, it is easy to see by definition that $\X\subseteq\bigcup_{n\ge 0}\res^n\X\subseteq\res\X$.
It remains to show that $\bigcup_{n\ge 0}\res^n\X$ is a resolving subcategory of $\mod R$.
But this is also easy to check.
\end{rem}

From its definition, we might think that there are not so many nonisomorphic indecomposable modules in $\res^n\X$.
But, the following two examples say that this guess is not right.

\begin{ex}
Let us consider the $1$-dimensional complete local hypersurface $R=\C[[x,y]]/(x^2)$ over the complex number field.
Then the subcategory $\res^1(xR)$ coincides with $\CM(R)$, and there exist infinitely many nonisomorphic indecomposable $R$-modules in $\res^1(xR)$.

Indeed, set
$$
I_n=
\begin{cases}
R & (n=0), \\
(x,y^n)R & (0<n<\infty), \\
xR & (n=\infty).
\end{cases}
$$
It follows from \cite[Example 6.5]{Y} that the set $\{ I_n\}_{0\le n\le\infty}$ consists of all the nonisomorphic indecomposable maximal Cohen-Macaulay $R$-modules.
For each integer $n$ with $0<n<\infty$, we have isomorphisms
$$
((x,y^n)R)/xR \overset{f}{\gets} R/xR \overset{g}{\to} xR,
$$
where $f$ sends the residue class of $a\in R$ in $R/xR$ to the residue class of $y^na$ in $((x,y^n)R)/xR$, and $g$ sends the residue class of $a\in R$ in $R/xR$ to $xa\in xR$.
Hence we see that there is an exact sequence
$$
0 \to xR \to I_n \to xR \to 0
$$
for $0\le n<\infty$, which implies that $I_n$ is in $\res^1(xR)$ for $0\le n\le \infty$.
On the other hand, $xR$ is a maximal Cohen-Macaulay $R$-module, and $\CM(R)$ is a resolving subcategory of $\mod R$ by Example \ref{ex}(\ref{cm}).
Therefore $\CM(R)$ coincides with $\res^1(xR)$.

The localization $R_\p$ at the prime ideal $\p=xR$ is singular because it has a (nonzero) nilpotent $x$.
Hence the local ring $R$ is not an isolated singularity, and thus there exist infinitely many isomorphism classes of indecomposable maximal Cohen-Macaulay $R$-modules by \cite[\S 10]{Au} or \cite[Corollary 2]{HL}.
\end{ex}

For a subcategory $\X$ of $\mod R$, we denote by $\ind\X$ (or $\ind_R\X$) the set of nonisomorphic indecomposable $R$-modules in $\X$.

\begin{ex}
Let $k$ be a field.
We consider the $2$-dimensional hypersurface $R=k[[x,y,z]]/(x^2)$.
Put $\p(f) =(x,y-zf)R$ for an element $f\in k[[z]]\subseteq R$.
Then $\p(f)$ is an indecomposable $R$-module in $\res^1(xR)$, and there exist uncountably many nonisomorphic indecomposable $R$-modules in $\res^1(xR)$.

Indeed, note that $xR$ is a maximal Cohen-Macaulay $R$-module.
Hence the $R$-regular element $y-zf$ is also $xR$-regular.
Note also that $\p(f)$ is isomorphic to $\Omega (xR/(y-zf)xR)$.
We can make the following pullback diagram:
$$
\begin{CD}
@. @. 0 @. 0 \\
@. @. @VVV @VVV \\
@. @. \p(f) @= \p(f) \\
@. @. @VVV @VVV \\
0 @>>> xR @>>> E @>>> R @>>> 0 \\
@. @| @VVV @VVV \\
0 @>>> xR @>{y-zf}>> xR @>>> xR/(y-zf)xR @>>> 0 \\
@. @. @VVV @VVV \\
@. @. 0 @. 0
\end{CD}
$$
Since the middle row splits, $E$ is isomorphic to $xR\oplus R$.
We get an exact sequence
$$
0 \to \p(f) \to xR\oplus R \to xR \to 0.
$$
The $R$-modules $xR$ and $xR\oplus R$ belong to $\res^0(xR)$, hence $\p(f)$ belongs to $\res^1(xR)$.
A similar argument to the proof of Claim 1 in \cite[Example 4.3]{count} shows that $\p(f)$ is an indecomposable $R$-module.
Thus, we obtain a map from $k[[z]]$ to $\ind(\res^1(xR))$ which is given by $f\mapsto\p(f)$.
Along the same lines as in the proofs of Claims 2 and 3 in \cite[Example 4.3]{count}, we can prove that this map is injective.
Since the set $k[[z]]$ is uncountably infinite, the assertion follows.
\end{ex}

For a subcategory $\X$ of $\mod R$ and a multiplicatively closed subset $S$ of $R$, we denote by $\X_S$ the subcategory of $\mod R_S$ consisting of all $R_S$-modules $X_S$ with $X\in\X$.
Our inductive construction of a resolving closure yields a relationship between a resolving closure and localization.

\begin{prop}\label{local}
Let $\X$ be a subcategory of $\mod R$, and let $S$ be a multiplicatively closed subset of $R$.
Then $(\res_R\X)_S$ is contained in $\res_{R_S}\X_S$.
\end{prop}

\begin{proof}
It is enough to show that $(\res_R^n\X)_S$ is contained in $\res_{R_S}^n\X_S$ for each integer $n\ge 0$.
We use induction on $n$.
The statement obviously holds when $n=0$.
Let $n\ge 1$, and take an $R$-module $M$ in $\res_R^n\X$.
Then there are a finite number of $R$-modules $M_1,\dots,M_t$ such that $M$ is a direct summand of $M_1\oplus\cdots\oplus M_t$ and that for each $1\le i\le t$ there exists an exact sequence of either of the following two forms:
\begin{align*}
& 0 \to A_i \to M_i \to B_i \to 0, \\
& 0 \to M_i \to A_i \to B_i \to 0
\end{align*}
where $A_i$ and $B_i$ are in $\res_R^{n-1}\X$.
Hence for each $1\le i\le t$ there is an exact sequence of either of the following two forms:
\begin{align*}
& 0 \to (A_i)_S \to (M_i)_S \to (B_i)_S \to 0, \\
& 0 \to (M_i)_S \to (A_i)_S \to (B_i)_S \to 0
\end{align*}
where $(A_i)_S$ and $(B_i)_S$ are in $(\res_R^{n-1}\X)_S$.
Induction hypothesis implies that $(\res_R^{n-1}\X)_S$ is contained in $\res_{R_S}^{n-1}\X_S$.
Since $M_S$ is a direct summand of $(M_1)_S\oplus\cdots\oplus(M_t)_S$, the $R_S$-module $M_S$ belongs to $\res_{R_S}^n\X_S$.
\end{proof}

Making use of the above result, we see that the nonfree locus of a subcategory is stable under taking its resolving closure.

\begin{cor}\label{nfres}
The equalities
$$
\NF(\res\X)=\NF(\add\X)=\NF(\X)
$$
hold for each subcategory $\X$ of $\mod R$.
\end{cor}

\begin{proof}
Note that there are inclusions $\res_R\X\supseteq\add_R\X\supseteq\X$ of subcategories of $\mod R$.
From this we see that there are inclusions $\NF(\res_R\X)\supseteq\NF(\add_R\X)\supseteq\NF(\X)$ of subsets of $\Spec R$.
We have only to show that $\NF(\res_R\X)$ is contained in $\NF(\X)$.

Let $\p$ be a prime ideal in $\NF(\res_R\X)$.
Then there is an $R$-module $Y\in\res_R\X$ such that $\p$ is in $\NF(Y)$.
The localization $Y_\p$ belongs to $(\res_R\X)_\p$, and to $\res_{R_\p}\X_\p$ by Proposition \ref{local}.
Assume that $\p$ is not in $\NF(\X)$.
Then for every $X\in\X$ the $R_\p$-module $X_\p$ is free.
Hence the subcategory $\X_\p$ of $\mod R_\p$ consists of all free $R_\p$-modules, and in particular, $\X_\p$ is resolving by Example \ref{ex}(\ref{proj}).
Therefore we have $\res_{R_\p}\X_\p=\X_\p$, and thus $Y_\p$ is a free $R_\p$-module.
But this contradicts the choice of $\p$.
Consequently, the prime ideal $\p$ must be in $\NF(\X)$, which completes the proof of the corollary.
\end{proof}

Using the above corollary, we can show that the nonfree locus of a subcategory is determined by the isomorphism classes of indecomposable modules in its resolving closure.

\begin{cor}\label{refment}
Let $\X$ be a subcategory of $\mod R$.
Then one has $\NF(\X)=\bigcup_{Y,Z\in\ind(\res\X)}\Supp\Ext_R^1(Y,Z)$.
\end{cor}

\begin{proof}
By Corollary \ref{nfres}, replacing $\X$ with $\res\X$, we may assume that the subcategory $\X$ is resolving.
Under this assumption, we have only to show the equality $\NF(\X)=\bigcup_{Y,Z\in\ind\X}\Supp\Ext_R^1(Y,Z)$.
If a prime ideal $\p$ is such that $\Ext_{R_\p}^1(Y_\p,Z_\p)\ne 0$ for some modules $Y,Z\in\X$, then $Y_\p$ is nonfree as an $R_\p$-module, hence $\p$ is in $\NF(\X)$.
Conversely, let $\p$ be a prime ideal in $\NF(X)$ for some $X\in\X$.
Then it follows from Proposition \ref{nfsupp} that $\Ext_{R_\p}^1(X_\p,\Omega X_\p)$ is nonzero.
Hence there are indecomposable summands $Y$ and $Z$ of $X$ and $\Omega X$ respectively such that $\Ext_{R_\p}^1(Y_\p,Z_\p)$ is nonzero.
The modules $Y,Z$ are in $\ind\X$.
\end{proof}

\section{Closed subsets of nonfree loci}

In this section, we study the structure of the nonfree locus of an $R$-module.
The main result of this section is concerning closed subsets of a nonfree locus (in the relative topology induced by the Zariski topology of $\Spec R$), which will often be referred in later sections.
We begin with the following lemma, which is proved by taking advantage of an idea used in the proof of \cite[Theorem 1]{HL}.

\begin{lem}\label{hl}
Let $R$ be a local ring with maximal ideal $\m$.
Let
$$
\sigma: 0\to L\overset{f}{\to} M\to N\to 0
$$
be an exact sequence of $R$-modules.
Let $x$ be an element in $\m$.
Then there is an exact sequence
$$
0\to L\overset{\binom{x}{f}}{\to} L\oplus M\to K\to 0.
$$
If this splits, then so does $\sigma$.
\end{lem}

\begin{proof}
There exists a homomorphism $(g,h):L\oplus M\to L$ such that $1=(g,h)\binom{x}{f}=xg+hf$.
Applying $\Hom_R(N,-)$ to $\sigma$, we have an exact sequence
$$
\Hom_R(N,N) \overset{\eta}{\to} \Ext_R^1(N,L) \overset{\Ext_R^1(N,f)}{\longrightarrow} \Ext_R^1(N,M),
$$
and get $\Ext_R^1(N,f)(\sigma)=(\Ext_R^1(N,f)\cdot\eta)(1)=0$.
Set $\xi=\Ext_R^1(N,g)$.
There are equalities $1=\Ext_R^1(N,xg+hf)=x\xi+\Ext_R^1(N,h)\cdot\Ext_R^1(N,f)$, so we obtain $\sigma=x\xi(\sigma)+\Ext_R^1(N,h)(\Ext_R^1(N,f)(\sigma))=x\xi(\sigma)$.
Hence $\sigma=x^i\xi^i(\sigma)$ for any $i\ge 1$, and therefore $\sigma\in\bigcap_{i\ge 1}\m^i\Ext_R^1(N,L)=0$ by virtue of Krull's intersection theorem.
Thus the exact sequence $\sigma$ splits.
\end{proof}

Using the above lemma, we prove the following proposition, which will play an essential role in the proofs of our main results.

\begin{prop}\label{const}
Let $X$ be an $R$-module.
Let $\p$ be a prime ideal in $\NF(X)$ and $x$ an element in $\p$.
Then there is a commutative diagram
\begin{equation}\label{cd}
\begin{CD}
\phantom{x}\sigma\ @. :\ @. 0 @>>> \Omega X @>{f}>> R^n @>>> X @>>> 0 \\
@. @. @. @V{x}VV @VVV @| \\
x\sigma\ @. :\ @. 0 @>>> \Omega X @>>> X_1 @>>> X @>>> 0
\end{CD}
\end{equation}
of $R$-modules with exact rows, and the following statements hold:
\begin{enumerate}[\rm (1)]
\item
$X_1\in\res^2X$,
\item
$V(\p)\subseteq\NF(X_1)\subseteq\NF(X)$,
\item
$D(x)\cap\NF(X_1)=\emptyset$.
\end{enumerate}
\end{prop}

\begin{proof}
Taking a free cover of $X$, we get an exact sequence
$$
\sigma: 0 \to \Omega X \overset{f}{\to} R^n \to X \to 0
$$
of $R$-modules.
Making a pushout diagram of $f:\Omega X\to R^n$ and the multiplication map $x:\Omega X\to\Omega X$, we obtain a commutative diagram \eqref{cd}.

(1) From the first row in \eqref{cd} we see that $\Omega X$ is in $\res^1X$.
It follows from the second row that $X_1$ is in $\res^2X$.

(2) Assume that $\p$ is not in $\NF(X_1)$.
Then $(X_1)_\p$ is free as an $R_\p$-module, and the exact sequence
$$
0 \to \Omega X_\p \overset{\binom{x}{f_\p}}{\to} \Omega X_\p\oplus R_\p^n \to (X_1)_\p \to 0
$$
splits.
Lemma \ref{hl} implies that $\sigma_\p$ is a split exact sequence, and so the $R_\p$-module $X_\p$ is free.
This contradicts the assumption that $\p$ is in $\NF(X)$.
Therefore $\p$ is in $\NF(X_1)$, and the set $V(\p)$ is contained in $\NF(X_1)$ by Corollary \ref{cl-spcl}(2).

Take a prime ideal $\q\in\NF(X_1)$.
Suppose that $\q$ is not in $\NF(X)$.
Then $X_\q$ is a free $R_\q$-module, and the localized exact sequences
\begin{align*}
\sigma_\q & :\  0 \to \Omega X_\q \to R^n_\q \to X_\q \to 0,\\
x\sigma_\q & :\  0 \to \Omega X_\q \to (X_1)_\q \to X_\q \to 0
\end{align*}
both split.
This implies that $(X_1)_\q$ is a free $R_\q$-module, which contradicts the choice of $\q$.
Thus $\NF(X_1)$ is contained in $\NF(X)$.

(3) Assume that the set $D(x)\cap\NF(X_1)$ is nonempty, and take a prime ideal $\q$ in $D(x)\cap\NF(X_1)$.
Then the element $x$ can be regarded as a unit of the local ring $R_\q$.
Localizing the diagram \eqref{cd} at $\q$, we see from the five lemma that $(X_1)_\q$ is a free $R_\q$-module.
Hence $\q$ is not in $\NF(X_1)$, which is a contradiction.
\end{proof}

Now we can prove one of the main results of this paper.

\begin{thm}\label{thm1}
For any $R$-module $X$ and any subset $W$ of $\NF(X)$ which is closed in $\Spec R$, there exists an $R$-module $Y\in\res X$ such that $W=\NF(Y)$.
\end{thm}

\begin{proof}
First of all, if $W$ is an empty set, then we can take $Y:=R$.
So, suppose that $W$ is nonempty.
Take an irreducible decomposition $W=V(\p_1)\cup\cdots\cup V(\p_n)$ of $W$.
Suppose that for each $1\le i\le n$ we can find an $R$-module $Y_i\in\res X$ such that $\NF(Y_i)$ coincides with $V(\p_i)$.
Then, putting $Y=Y_1\oplus\cdots\oplus Y_n$, which belongs to $\res X$, we easily see that $W$ is equal to $\NF(Y)$.
So, we can assume without loss of generality that $W$ is an irreducible closed subset of $\Spec R$; we write $W=V(\p)$ for some $\p\in\Spec R$.

If $V(\p)$ coincides with $\NF(X)$, then we can take $Y:=X$.
So assume that $V(\p)$ is strictly contained in $\NF(X)$.
Then there is a prime ideal $\q\in\NF(X)$ which is not in $V(\p)$.
Hence there exists an element $x\in\p$ which is not in $\q$.
For this element $x$ of $R$, let $X_1$ be an $R$-module satisfying the three conditions in Proposition \ref{const}.
Then it is obvious that $X_1$ is in $\res X$.
Since $\q$ is in $D(x)$, it does not belong to $\NF(X_1)$.
Thus we have $V(\p)\subseteq\NF(X_1)\subsetneq\NF(X)$.
If $V(\p)$ coincides with $\NF(X_1)$, then we can take $Y:=X_1$.
So we assume that $V(\p)$ is strictly contained in $\NF(X_1)$.
Then, a similar argument to the above shows that there exists an $R$-module $X_2\in\res X$ which satisfies $V(\p)\subseteq\NF(X_2)\subsetneq\NF(X_1)\subsetneq\NF(X)$.

According to Corollary \ref{cl-spcl}(1), all nonfree loci are closed subsets of $\Spec R$.
Since $\Spec R$ is a noetherian space, every descending chain of closed subsets stablizes.
This means that the above procedure to construct modules $X_i$ cannot be iterate infinitely many times.
Hence there exists an $R$-module $Y\in\res X$ such that $V(\p)$ coincides with $\NF(Y)$.
\end{proof}

A very special case of this theorem has already been obtained by the author; see \cite[Lemma 3.4]{count}.

From the above theorem we see that the nonfree locus of a given nonfree module has an irreducible decomposition by the nonfree loci of a finite number of modules in its resolving closure.

\begin{cor}
For every nonfree $R$-module $X$ there exists a decomposition
$$
\NF(X)=\NF(Y_1)\cup\cdots\cup\NF(Y_n)
$$
with $Y_1,\dots,Y_n\in\res X$ such that $\NF(Y_1),\dots,\NF(Y_n)$ are irreducible closed subsets of $\Spec R$.
\end{cor}

\begin{proof}
Since $\NF(X)$ is a nonempty closed subset of $\Spec R$, we have a decomposition
$$
\NF(X)=V(\p_1)\cup\cdots\cup V(\p_n)
$$
for some prime ideals $\p_1,\dots,\p_n$.
We apply Theorem \ref{thm1} to each $V(\p_i)$ to see that there is an $R$-module $Y_i\in\res X$ such that $V(\p_i)$ coincides with $\NF(Y_i)$.
Then each $\NF(Y_i)$ is irreducible and we have $\NF(X)=\NF(Y_1)\cup\cdots\cup\NF(Y_n)$.
\end{proof}

We might think that the above corollary predicts that all closed subsets of $\Spec R$ are the nonfree loci of some modules.
But, as the proposition below says, this statement does not hold.
Here, for a subset $W$ of $\Spec R$, we denote by $\min W$ the set of minimal elements of $W$ with respect to inclusion relation.

\begin{prop}
Let $W$ be a nonempty closed subset of $\Spec R$.
Then the following are equivalent:
\begin{enumerate}[\rm (1)]
\item
One has $W=\NF(X)$ for some $R$-module $X$;
\item
One has $W=\NF(X)$ for some $R$-module $X\in\res_R(\bigoplus_{\p\in\min W}R/\p)$;
\item
For every $\p\in W$, the local ring $R_\p$ is not a field.
\end{enumerate}
\end{prop}

\begin{proof}
(2) $\Rightarrow$ (1): This implication is trivial.

(1) $\Rightarrow$ (3): Let $\p$ be a prime ideal in $W$.
Then $X_\p$ is not a free $R_\p$-module.
In particular, $R_\p$ is not a field.

(3) $\Rightarrow$ (2): Take an irreducible decomposition $W=V(\p_1)\cup\cdots\cup V(\p_n)$ of $W$.
Fix an integer $i$ with $1\le i\le n$.
By assumption, the local ring $R_{\p_i}$ is not a field.
It is easy to see that $\p_i$ belongs to $\NF_R(R/\p_i)$, hence $V(\p_i)$ is contained in $\NF_R(R/\p_i)$ by Corollary \ref{cl-spcl}(2).
Theorem \ref{thm1} implies that there exists an $R$-module $Y_i\in\res_R(R/\p_i)$ such that $V(\p_i)=\NF(Y_i)$.
Setting $Y=Y_1\oplus\cdots\oplus Y_n$, we see that $Y$ is in $\res_R(\bigoplus_{i=1}^nR/\p_i)$ and that $W$ coincides with $\NF(Y)$.
\end{proof}

\section{Walks in resolving subcategories}

In this section, we investigate the structure of the resolving closure of an $R$-module by means of the inductive construction of the resolving closure which we obtained in Section 3.
More precisely, let $X$ be an $R$-module.
For an $R$-module $Y\in\res X$, we consider how many resolving operations are needed to take to construct $Y$ from $X$.
Here, resolving operations mean extensions and kernels of epimorphisms.
For this purpose, we introduce the following invariant which measures the minimum number of required resolving operations.
This is an imitation of a level in a triangulated category defined in \cite{ABIM}.

\begin{dfn}\label{deflevel}
For two $R$-modules $X$ and $Y$, we define
$$
\level(X,Y)=\inf\{ n\ge 0\mid Y\in\res_R^nX\}.
$$
\end{dfn}

\begin{rem}
Let $X$ be an $R$-module.
\begin{enumerate}[(1)]
\item
One has $\level(X,Y)=0$ for every $R$-module $Y\in\res^0X=\add(X\oplus R)$.
In particular, $\level(X,X)=\level(X,R)=\level(X,0)=0$.
\item
One has $\level(X,Y)<\infty$ for an $R$-module $Y$ if and only if $Y$ belongs to $\res X$.
\end{enumerate}
\end{rem}

In general, the invariant $\level(-,-)$ does not induce a distance function.
However, it satisfies the triangle inequality.

\begin{prop}\label{3kaku}
Let $X,Y,Z$ be $R$-modules.
\begin{enumerate}[\rm (1)]
\item
Let $m,n$ be nonnegative integers.
If $Y\in\res^mX$ and $Z\in\res^nY$, then $Z\in\res^{m+n}X$.
\item
The inequality $\level(X,Z)\le\level(X,Y)+\level(Y,Z)$ holds.
\end{enumerate}
\end{prop}

\begin{proof}
(1) Let us prove this assertion by induction on $n$.

When $n=0$, the module $Z$ is in $\add(Y\oplus R)$.
Note that both $Y$ and $R$ are in $\res^mX$.
Since $\res^mX$ is an additive closure, it contains $\add(Y\oplus R)$.
Hence $Z$ belongs to $\res^mX=\res^{m+n}X$.

Let $n\ge 1$.
By definition, there are a finite number of $R$-modules $M_1,\dots,M_s$ such that $Z$ is a direct summand of $M_1\oplus\cdots\oplus M_s$ and that for each $1\le i\le s$ there exists an exact sequence of either of the following two forms:
\begin{align*}
& 0 \to A_i \to M_i \to B_i \to 0, \\
& 0 \to M_i \to A_i \to B_i \to 0
\end{align*}
where $A_i,B_i\in\res^{n-1}Y$.
Induction hypothesis implies that the modules $A_i,B_i$ are in $\res^{m+n-1}X$ for $1\le i\le s$.
Hence each $M_i$ is in $\res^{m+n}X$, and therefore so is $Z$.

(2) Set $p=\level(X,Y)$ and $q=\level(Y,Z)$.
Then $Y$ is in $\res^pX$ and $Z$ is in $\res^qY$.
The assertion (1) implies that $Z$ is in $\res^{p+q}X$, which says that $\level(X,Z)\le p+q=\level(X,Y)+\level(Y,Z)$.
\end{proof}

Let $Z$ be a subset of $\Spec R$.
For a prime ideal $\p$ in $Z$, we define the {\em height} of $\p$ with respect to $Z$ as the supremum of $\height(\p/\q)$ where $\q$ runs through all prime ideals in $Z$ that are contained in $\p$.
We denote it by $\height_Z(\p)$.

\begin{rem}
The following statements are straightforward.
\begin{enumerate}[(1)]
\item
One has $\height_{\Spec R}(\p)=\height\p$ for any $\p\in\Spec R$.
\item
For a prime ideal $\p$ in a subset $Z$ of $\Spec R$, it holds that $0\le\height_Z(\p)\le\height\p$.
\item
Let $Z$ be a closed subset of $\Spec R$, and let $I$ be the defining ideal of $Z$.
For each prime ideal $\p\in Z$, one has $\height_Z(\p)=\height(\p/I)$.
\item
Let $\p,\q$ be prime ideals in a subset $Z$ of $\Spec R$.
If $\p\subseteq\q$, then $\height_Z(\p)\le\height_Z(\q)$.
\item
Let $Z,W$ be subsets of $\Spec R$.
If $Z\subseteq W$, then $\height_Z(\p)\le\height_W(\p)$ for any $\p\in Z$.
\item
Let $R$ be a local ring with maximal ideal $\m$.
Then the equality $\height_Z(\m)=\dim Z$ holds for every subset $Z$ of $\Spec R$ containing $\m$.
(Recall that the {\em dimension} $\dim Z$ of a subset $Z$ of $\Spec R$ is defined as the supremum of $\dim R/\p$ where $\p$ runs over all prime ideals in $Z$.)
\item
Let $Z$ be a subset of $\Spec R$ and let $\p$ be a prime ideal in $Z$.
Then $\height_Z(\p)=0$ if and only if $\p$ is minimal in $Z$.
\end{enumerate}
\end{rem}

Now we state and prove one of the main results of this paper.

\begin{thm}\label{thm2}
Let $X$ be an $R$-module and let $\p$ be a prime ideal in $\NF(X)$.
Then there exists an $R$-module $Y\in\res X$ satisfying the following three conditions:
\begin{enumerate}[\rm (1)]
\item
$\level(X,Y)\le 2\height_{\NF(X)}(\p)$,
\item
$\p\in\NF(Y)$,
\item
$\height_{\NF(Y)}(\p)=0$.
\end{enumerate}
\end{thm}

\begin{proof}
We prove the theorem by induction on $n:=\height_{\NF(X)}(\p)$.
(Note that $\height_{\NF(X)}(\p)$ is finite because $R$ is a noetherian ring.)

When $n=0$, we set $Y:=X$.
Then $Y$ is in $\res^0X$, so we have $\level(X,Y)=0=2\height_{\NF(X)}(\p)$.
We also have $\p\in\NF(Y)$ and $\height_{\NF(Y)}(\p)=0$.

When $n\ge 1$, put
$$
S=\{\q\in\NF(X)\mid\height_{\NF(X)}(\q)=0\}.
$$
Corollary \ref{cl-spcl}(1) implies that $\NF(X)$ is a closed subset of $\Spec R$.
Letting $I$ be the defining ideal of $\NF(X)$, we have $\height_{\NF(X)}(\q)=\height(\q/I)$ for every $\q\in\NF(X)$.
Hence $S$ coincides with the set of minimal prime ideals of $I$, and therefore $S$ is a finite set.
As $n$ is positive now, the prime ideal $\p$ is not contained in all prime ideals in $S$.
By prime avoidance, we can choose an element $x\in\p$ which is not contained in all prime ideals in $S$.

For this element $x$, take an $R$-module $X_1$ which satisfies the conditions in Proposition \ref{const}.
Namely, the module $X_1$ satisfies the following three conditions:
\begin{align*}
& X_1\in\res^2X,\\
& V(\p)\subseteq\NF(X_1)\subseteq\NF(X),\\
& D(x)\cap\NF(X_1)=\emptyset.
\end{align*}
Hence $\level(X,X_1)\le 2$ and $X_1\in\res X$.
Since $S$ is contained in $D(x)$, we have $S\cap\NF(X_1)=\emptyset$.

Let $\q$ be a prime ideal in $\NF(X_1)$ which is contained in $\p$.
Then $\q$ does not belong to $S$, so $\height_{\NF(X)}(\q)>0$.
Hence $\height(\q/\r)>0$ for some prime ideal $\r\in\NF(X)$ which is contained in $\q$.
There are inequalities
$$
\height(\p/\q)<\height(\p/\q)+\height(\q/\r)\le\height(\p/\r)\le\height_{\NF(X)}(\p)=n.
$$
Therefore we have $\height_{\NF(X_1)}(\p)<n$.
The induction hypothesis implies that there exists an $R$-module $Y\in\res X_1$ such that $\level(X_1,Y)\le 2\height_{\NF(X_1)}(\p)$, that $\p\in\NF(Y)$ and that $\height_{\NF(Y)}(\p)=0$.
According to Proposition \ref{3kaku}(2), there are inequalities
\begin{align*}
\level(X,Y) & \le\level(X,X_1)+\level(X_1,Y) \\
& \le 2+2\height_{\NF(X_1)}(\p) \\
& \le 2+2(n-1)=2n=2\height_{\NF(X)}(\p).
\end{align*}
Thus the proof of the theorem is completed.
\end{proof}

Applying the above theorem to a local ring $R$, we get the following result.
This result contains Theorem A from the introduction.

\begin{cor}\label{locthm2}
Let $R$ be a local ring.
Then for every nonfree $R$-module $X$, there exists a nonfree $R$-module $Y$ in $\res X$ satisfying the following conditions:
\begin{enumerate}[\rm (1)]
\item
$\level(X,Y)\le 2\dim\NF(X)$,
\item
$Y$ is free on the punctured spectrum of $R$.
\end{enumerate}
\end{cor}

\begin{proof}
Let $\m$ be the unique maximal ideal of $R$.
We observe that $\m$ is in $\NF(X)$.
Letting $\p=\m$ in Theorem \ref{thm2}, we see that there is an $R$-module $Y\in\res X$ such that $\level(X,Y)\le 2\height_{\NF(X)}(\m)$, that $\m\in\NF(Y)$ and that $\height_{\NF(Y)}(\m)=0$.
These three conditions imply that the inequality $\level(X,Y)\le 2\dim\NF(X)$ holds, that $Y$ is a nonfree $R$-module and that $Y_\p$ is a free $R_\p$-module for every prime ideal $\p\ne \m$, respectively.
\end{proof}

Restricting the above corollary to the Cohen-Macaulay case, we obtain the following result on maximal Cohen-Macaulay modules.

\begin{cor}
Let $R$ be a Cohen-Macaulay local ring.
Then for any nonfree maximal Cohen-Macaulay $R$-module $X$, there exists a nonfree maximal Cohen-Macaulay $R$-module $Y$ satisfying the following two conditions:
\begin{enumerate}[\rm (1)]
\item
$\level(X,Y)\le 2\dim\Sing R$,
\item
$Y$ is free on the punctured spectrum of $R$.
\end{enumerate}
\end{cor}

\begin{proof}
By virtue of Corollary \ref{locthm2}, we find an $R$-module $Y\in\res X$ which is free on the punctured spectrum of $R$ and satisfies the inequality $\level(X,Y)\le 2\dim\NF(X)$.
Since $X$ is in $\CM(R)$ and $\CM(R)$ is resolving by Example \ref{ex}(\ref{cm}), the module $Y$ is also in $\CM(R)$, that is, $Y$ is maximal Cohen-Macaulay.
Since $\NF(\CM(R))=\Sing R$ by Example \ref{cmex}, the assertion follows.
\end{proof}

Forgetting the first condition on the module $Y$ in Corollary \ref{locthm2}, we obtain the following result.

\begin{cor}\label{existnf}
Let $R$ be a local ring and $\X$ a resolving subcategory of $\mod R$.
If there exists a nonfree $R$-module in $\X$, then there exists a nonfree $R$-module in $\X$ which is free on the punctured spectrum of $R$.
\end{cor}

\begin{rem}
In the case where $R$ is a Cohen-Macaulay local ring, a nonfree $R$-module in $\X:=\CM(R)$ which is free on the punctured spectrum can be constructed explicitly as follows.
Let $R$ be a $d$-dimensional Cohen-Macaulay local ring with maximal ideal $\m$.
Then it is well-known and easy to see that there exists a nonfree $R$-module in $\CM(R)$ if and only if $R$ is singular.
When this is the case, the $R$-module $\Omega^d(R/\m)$ is a nonfree $R$-module in $\CM(R)$ which is free on the punctured spectrum of $R$.
\end{rem}

Applying Corollary \ref{existnf} to the resolving subcategory $\X=\G(R)$ (cf. Example \ref{ex}(\ref{tr})), we have the following.

\begin{cor}
Let $R$ be a local ring.
If there exists a nonfree totally reflexive $R$-module, then there exists a nonfree totally reflexive $R$-module which is free on the punctured spectrum.
\end{cor}

\begin{rem}
A local ring over which all totally reflexive modules are free is called {\em G-regular}.
G-regular local rings have been studied by several authors.
One of the main problems for G-regular local rings is to establish necessary and/or sufficient conditions for a given local ring to be G-regular.
For the details of G-regular local rings, see \cite{greg}.
The above corollary should give some contribution to this problem.
\end{rem}

\section{Resolving subcategories of countable type}

In this section, we investigate resolving subcategories in which there exist only countably many nonisomorphic indecomposable modules.
The following proposition plays a key role for this goal, which is proved by using Theorem \ref{thm1}.
It actually gives a refinement of one inclusion in the equality given in Corollary \ref{refment}.

\begin{prop}\label{inc}
For a subcategory $\X$ of $\mod R$ one has an inclusion of sets:
$$
\NF(\X)\subseteq\left\{\sqrt{\Ann\Ext_R^1(Y,Z)}\ \Bigg |\ Y,Z\in\ind(\res\X)\right\}.
$$
\end{prop}

\begin{proof}
Let $\p$ be a prime ideal in $\NF(\X)$.
Then $\p$ is in $\NF(X)$ for some $R$-module $X\in\X$.
As $\NF(X)$ is specialization-closed by Corollary \ref{cl-spcl}(2), the irreducible set $V(\p)$ is contained in $\NF(X)$.
According to Theorem \ref{thm1}, there exists an $R$-module $Y\in\res X$ such that $V(\p)$ coincides with $\NF(Y)$.
Since $\NF(Y)=\Supp\Ext_R^1(Y,\Omega Y)=V(\Ann\Ext_R^1(Y,\Omega Y))$ by Proposition \ref{nfsupp}, the prime ideal $\p$ is equal to $\sqrt{\Ann\Ext_R^1(Y,\Omega Y)}$.
Take indecomposable decompositions $Y=\bigoplus_{i=1}^mY_i$ and $\Omega Y=\bigoplus_{j=1}^nZ_j$.
Then we have
\begin{align*}
\p & =\sqrt{\Ann\Ext_R^1(Y,\Omega Y)} \\
& =\sqrt{\Ann\Ext_R^1(\bigoplus_{i=1}^mY_i,\bigoplus_{j=1}^nZ_j)} \\
& =\bigcap_{1\le i\le m,1\le j\le n}\sqrt{\Ann\Ext_R^1(Y_i,Z_j)}.
\end{align*}
Since $\p$ is a prime ideal, $\p$ is equal to $\sqrt{\Ann\Ext_R^1(Y_a,Z_b)}$ for some integers $a,b$.
As $\res X$ is a resolving subcategory of $\mod R$, the $R$-modules $Y_a,Z_b$ are in $\res X$, hence in $\res\X$ and therefore in $\ind(\res\X)$.
Thus we obtain the desired inclusion.
\end{proof}

A very special case of the above proposition has already been obtained by the author; see \cite[Proposition 3.5]{count}.

\begin{dfn}
We say that a subcategory $\X$ of $\mod R$ has {\em coubtable type} if the set $\ind\X$ is countable.
\end{dfn}

We say that a Cohen-Macaulay local ring $R$ has {\em countable Cohen-Macaulay representation type} if $\CM(R)$ has countable type.

The result below is a direct consequence of Proposition \ref{inc}.

\begin{cor}\label{gyaku}
Let $\X$ be a subcategory of $\mod R$.
If $\res\X$ has countable type, then $\NF(\X)$ is at most a countable set.
\end{cor}

The converse of this corollary does not necessarily hold.
Indeed, we have the following example.

\begin{ex}
We consider a $1$-dimensional local hypersurface $R=\C[[x,y]]/(x^4+y^5)$.
Let $\m=(x,y)$ be the maximal ideal of $R$, and set $\X=\CM(R)$.
Then, since this ring $R$ is an integral domain of dimension $1$, we have $\NF(\X)=\Sing R=\{\m\}$ (cf. Example \ref{cmex}).
In particular, the set $\NF(\X)$ is finite, hence at most countable.
Since $\X$ is resolving by Example \ref{ex}(\ref{cm}), we have $\res\X=\X$.
This subcategory $\X$ does not have countable type by virtue of the classification theorem \cite[Theorem B]{BGS} of hypersurfaces of finite and countable Cohen-Macaulay representation type.
\end{ex}

The lemma below is proved by using so-called countable prime avoidance; see \cite[Lemma 2.2]{count} for the proof.

\begin{lem}\label{ctlem}
Let $R$ be a local ring with residue field $k$, and assume either that $R$ is complete or that $k$ is uncountable.
Let $Z$ be a specialization-closed subset of $\Spec R$.
If $Z$ is at most countable, then $\dim Z\le 1$.
\end{lem}

Corollaries \ref{gyaku}, \ref{cl-spcl}(2) and Lemma \ref{ctlem} yield the following theorem, which is one of the main results of this paper.

\begin{thm}\label{thm3}
Let $R$ be a local ring with residue field $k$, and assume either that $R$ is complete or that $k$ is uncountable.
Let $\X$ be a subcategory of $\mod R$ such that $\res\X$ has countable type.
Then $\dim\NF(\X)\le 1$.
\end{thm}

Combining this theorem with Corollary \ref{locthm2} gives the following result.

\begin{cor}
Let $R$ be a local ring with residue field $k$, and assume either that $R$ is complete or that $k$ is uncountable.
Let $X$ be a nonfree $R$-module such that $\res X$ has countable type.
Then there exists a nonfree $R$-module $Y\in\res X$ which is free on the punctured spectrum of $R$ and satisfies $\level(X,Y)\le 2$.
\end{cor}

\begin{rem}
Along the lines in the proof of Theorem \ref{thm2}, we can actually construct such a module $Y$ as in the above corollary.
Let $\m$ be the unique maximal ideal of $R$.
Since $X$ is a nonfree $R$-module, $\m$ belongs to $\NF(X)$.
Theorem \ref{thm3} guarantees that $\height_{\NF(X)}(\m)=\dim\NF(X)\le 1$.
If $\height_{\NF(X)}(\m)=0$, then $X$ is free on the punctured spectrum, so we can take $Y:=X$.
In this case we have $\level(X,Y)=0$.
If $\height_{\NF(X)}(\m)=1$, then the proof of Theorem \ref{thm2} implies that there exists an element $x\in\m$ which is not in each prime ideal $\q\in\NF(X)$ with $\height_{\NF(X)}(\q)=0$.
Applying Proposition \ref{const} to this element $x$, we obtain an $R$-module $X_1$ satisfying the three conditions in the proposition.
The proof of Theorem \ref{thm2} shows that $\height_{\NF(X_1)}(\m)<1$, which implies that $X_1$ is free on the punctured spectrum of $R$.
Thus we can take $Y:=X_1$.
In this case we have $\level(X,Y)\le 2$.
\end{rem}

We immediately get the following corollary from Theorem \ref{thm3}.
This result is nothing but Theorem B from the introduction.

\begin{cor}\label{corthm3}
Let $R$ be a local ring with residue field $k$, and assume either that $R$ is complete or that $k$ is uncountable.
Let $\X$ be a resolving subcategory of $\mod R$ of countable type.
Then $\dim\NF(\X)\le 1$.
\end{cor}

Applying this corollary to the subcategory of maximal Cohen-Macaulay modules over a Cohen-Macaulay local ring (cf. Example \ref{ex}(\ref{cm})), we can recover a theorem of Huneke and Leuschke.

\begin{cor}\cite[Theorem 1.3]{HL}\cite[Theorem 2.4]{count}
Let $R$ be a Cohen-Macaulay local ring of countable Cohen-Macaulay representation type.
Assume either that $R$ is complete or that the residue field is uncountable.
Then $\dim\Sing R\le 1$.
\end{cor}

Applying Corollary \ref{corthm3} to the subcategory of totally $C$-reflexive modules where $C$ is a semidualizing module (cf. Example \ref{ex}(\ref{tcr})), we obtain a refinement of the main theorem of \cite{count}.

\begin{cor}{\rm (cf. \cite[Theorem 3.6]{count})}
Let $R$ be a local ring which either is complete or has uncountable residue field.
Let $C$ be a semidualizing $R$-module.
Suppose that $\G_C(R)$ has countable type.
Then $\dim\NF(\G_C(R))\le 1$.
\end{cor}


\end{document}